\documentclass[12pt]{article}
\usepackage{amsmath}
\usepackage{amsthm}
\newtheorem{theorem}{Theorem}[section]
\newtheorem{lemma}[theorem]{Lemma}

\newtheorem{definition}[theorem]{Definition}

\newenvironment{pr}{\noindent{\bf Proof. }\rm }{\hfill{\tiny$\blacksquare$}}
\usepackage{amssymb}

\newcommand{\nd}{\mathbb{N}^d}
\newcommand{\mf}[1]{\mathbf{#1}}
\newcommand{\su}{ \sum \limits _{\mathbf{k} \leq \mathbf{n}}}
\newcommand{\ma}{ \max \limits _{\mathbf{k} \leq \mathbf{n}}}
\newcommand{\wt}[1]{\widetilde{#1}}
\newcommand{\x}{X_{\mathbf{k}}}
\newcommand{\sgn}{\operatorname{sgn}}
\begin{document}

\centerline{\large \bf AN APPROACH TO COMPLETE CONVERGENCE THEOREMS}

\smallskip
\centerline{\large \bf  FOR  DEPENDENT RANDOM FIELDS }

\smallskip
\centerline{\large \bf  VIA APPLICATION  OF FUK-NAGAEV INEQUALITY}

\bigskip
\centerline{\bf  Zbigniew A. Lagodowski}

\centerline{ Department of Mathematics, Lublin University of Technology,}

\centerline { Nadbystrzycka 38D, 20-618 Lublin, Poland}

\bigskip\bigskip
\noindent
{\bf Abstract:}
Let  $\{ X_{\bf n}, {\bf n}\in \mathbb{N}^d \}$ be a random field i.e. a family of random variables indexed by
$\mathbb{N}^d $, $d\ge 2$.  Complete convergence, convergence rates  for non identically distributed, negatively
dependent and martingale random
 fields are studied by application of Fuk-Nagaev inequality. The results are proved in asymmetric convergence case i.e.
 for the norming sequence equal $n_1^{\alpha_1}\cdot n_2^{\alpha_2}\cdot\ldots\cdot n_d^{\alpha_d}$,
 where $(n_1,n_2,\ldots, n_d)=\mathbf{n} \in  \mathbb{N}^d$ and $\min\limits_{1\leq i \leq d}\alpha_i \geq \frac{1}{2}.$



\bigskip\bigskip

\section{Introduction.}We observe again an interest in complete convergence theorems, which are discussed  for weighted sum of
dependent random variables, sums of random numbers of random variables, arrays of random variables or random fields. \\
We  will consider random variables on probability space $(\Omega, \mathfrak{F}, P)$, indexed by lattice points, i.e.
by index set
$\mathbb{N}^{d}, $ $d\geq 2$.
The elements of $\mathbb{N}^{d}$
denote: ${\bf m} = (m_1, m_2, \dots , m_d)$, ${\bf n} = (n_1, n_2, \dots , n_d)$ etc.,  we order them by coordinate wise
ordering: $ {\bf m} \le {\bf n} \quad \text{ iff} \quad m_i \le n_i ,  i = 1,2 \dots , d.$
We mean $ {\bf n} \rightarrow +\infty \quad  \text{iff} \quad \min\limits_{1\le i\le d}n_i \rightarrow +\infty.$
A family of random variables  $\{X_{\mathbf{n}}, \mathbf{n} \in \mathbb{N}^{d}\}$ we also call a random field, furthermore
denote
$S_{\mathbf{n}}=\sum \limits_{\mathbf{k} \leq \mathbf{n}}X_{\mathbf{k}}.$\\
This article is inspired by paper of Gut and Stadm\"{u}ller \cite{GutHsu}, where authors have studied Baum-Katz type theorems
 and   obtained  very general results for fields of independent identically distributed random variables while
 the normalizing sequence
 depends on different powers of different coordinates, i.e. they have studied convergence of the sums
\begin{equation}
\sum\limits_{\mathbf{n}}\mid \mathbf{n}\mid ^{\alpha_1r-2}
P(\max\limits_{\mathbf{k}\leq\mathbf{n}}\mid S_{\mathbf{k}}\mid >\mid \mathbf{n}^{\alpha}\mid \varepsilon),
\end{equation}
where $ \boldsymbol{\alpha} =(\alpha_1,\alpha_2,...,\alpha_d)\in (\frac{1}{2},1\rangle^d$, coordinates $\alpha_i$ are arranged
in non-decreasing order, $\alpha_1r \geq 1$ and $\mid \mathbf{n}^{\boldsymbol{\alpha}}\mid =
n_1^{\alpha_1}\cdot...\cdot n_d^{{\alpha}_d}$ or the case of convergence of
\begin{equation}
\sum\limits_{\mathbf{n}}\mid \mathbf{n}\mid ^{(r/2)-2}P(\max \limits_{\mathbf{k} \leq \mathbf{n}}\mid S_{\mathbf{n}}\mid \geq
\sqrt{ \prod\limits_{i=1}^pn_i\log(\prod\limits_{i=1}^pn_i)}\prod\limits_{i=p+1}^dn_i^{\alpha_i}\varepsilon),
\end{equation}
where $\alpha_1 =\frac{1}{2}$, $p=\max\{k:\alpha_k=\alpha_1\}$ and $r \geq 2$.\\
In the cited paper of Gut and Stadm\"{u}ller \cite{GutHsu} one can find the review and the comments of the so far
obtained  results and further references .
The crucial step in the proofs of the above mention theorems is based on symmetrization/desymmetrization and
Kahane-Hoffmann-J{\o}rgensen (K-H-J) inequality. K-H-J inequality is very sharp but strictly connected
 to independence of random variables.  In the proofs of Baum-Katz type theorems such a strong inequality is not
 needed, we can apply weaker one with an attribute of K-H-J inequality and at same time valid for dependent random variables.
  Fuk-Nagaev inequality is playing essential role in the proof of  such inequalities. Thus, by that approach, we are going to
  extend or give a compliments of some results of  Peligrad \cite{pelgut}, Gut et al.\cite{GutHsu}, \cite{Gut Asy},
  Kuczmaszewska et al. \cite {Kucz}.
Also, we will be able to extend the results of Ghosal et al. \cite{gho}, Sung \cite {Sung}, Dehua et al \cite{dehua}
to the random fields. Our result for martingale random field seems be a little bit more general even in one dimension case $(d=1)$
than the following result of Ghosal and Chandra (cf. Theorem 1(b) and Theorem 2 of \cite {gho})\\
\begin{theorem}
Let $\{(X_{nk}, \mathfrak{F}_{nk}), k \geq 1 \}$  be a sequence of square integrable martingale differences. Suppose, that
 there exist constants $(M_n)$ such that $\sum \limits _{k=1}^{\infty}E\left(X_{nk}^2\mid \mathfrak{F}_{n,k-1}\right) \leq
 M_n$ a.s., where $\mathfrak{F}_{n0}$ is trivial for all n. Let $(c_n)$ be a sequence of nonnegative numbers satisfying
  $\sum_{n=1}^{\infty}c_nM_n^\lambda <\infty $ for some $\lambda >0$ and
\begin{equation}\bigwedge \limits_{\varepsilon >0}\quad
\sum \limits_{n=1}^{\infty}c_n\sum\limits_{k=1}^{\infty}P\left(\mid X_{nk}\mid >\varepsilon \right)<\infty \quad,
\end{equation}
then
\begin{equation} \bigwedge \limits_{\varepsilon >0} \quad
\sum \limits_{n=1}^{\infty}c_n\sum\limits_{k=1}^{\infty}P\left(\sup\limits_{k>1}\mid \sum\limits _{i=1}^kX_{ni}\mid >
\varepsilon \right)<\infty \quad.
\end{equation}
\end{theorem}
\medskip

In order to formulate our main results we recall some definitions.
\begin{definition}
 A finite family of random variables   $\{X_{\mathbf{j}},\mathbf{1} \leq \mathbf{j} \leq \mathbf{n}\}$ is said to be
  negatively dependent  (ND) if
\begin{equation*}
P[\bigcap\limits_{\mathbf{j} \leq \mathbf{n}}(X_{\mathbf{j}} \leq x_{\mathbf{j}})] \leq \prod\limits_
{\mathbf{j} \leq \mathbf{n}}P(X_{\mathbf{j}} \leq x_{\mathbf{j}})
\end{equation*}
and
\begin{equation*}
P[\bigcap\limits_{\mathbf{j} \leq \mathbf{n}}(X_{\mathbf{j}} > x_{\mathbf{j}})] \leq \prod
\limits_{\mathbf{j} \leq \mathbf{n}}P(X_{\mathbf{j}} > x_{\mathbf{j}})
\end{equation*}
for  $x_\mathbf{i} \in \mathbb{R}$, $\mathbf{i} \leq \mathbf{n}$\\
\end {definition}
\noindent
An infinite family is ND if every finite subfamily is ND.\\
\begin{definition}
The family of random variables   $\{X_{\mathbf{j}}, \mathbf{j} \in \mathbb{N}^d\}$ is said to be
negatively associated (NA)  if

\begin{equation*}
cov(f(X_{\mathbf{j}}, \mathbf{j} \in S), g(X_{\mathbf{i}}, \mathbf{i} \in T)) \leq 0
\end{equation*}
for every disjoint subset $S,T\subset \mathbb{N}^d$, where  $f(X_{\mathbf{j}}, \mathbf{j} \in S) $ and
$g(X_{\mathbf{i}}, \mathbf{i} \in T) $
are coordinatwise increasing functions and the covariance exist.
\end{definition}
\noindent
An infinite family is NA if every finite subfamily is NA.\\
Since, we are going to prove  results for  non-identically distributed random variables, the following conditions allow us
 to formulate them in simple form as in i.i.d. case.
\medskip
\begin{definition} Random variables $\{ X_{\bf k},  {\bf k} \in \mathbb{N}^d\}$  are weakly mean  bounded (WMB)
 by random variable
$\xi$ (possibly defined on different probability space)
 iff there exist some constants $\kappa_1, \kappa_2 >0$, ${\bf n_0}\in \mathbb{N}^d $ and $x_0>0$ such that for every $x>x_0$ and
 ${\bf n}\geq {\bf n_0}$, ${\bf n} \in \mathbb{N}^d $
$$
\kappa_2\cdot P(\mid \xi\mid >x) \leq {1\over \mid {\bf n}\mid  }\sum_{{\bf k}\leq {\bf n}} P(\mid X_{\bf k}\mid >x )
\leq  \kappa_1\cdot P(\mid \xi\mid >x) $$
\end{definition}

\medskip
\noindent
 If only the right hand side inequality is satisfied, we say that the random field
$\{X_{\bf n}, {\bf n}\in \mathbb{N}^d\}$ and the random variable $\xi$ satisfy weak mean dominance (WMD) condition.
WMB condition seems to be  very natural one and not very restricted, e.g. regular cover condition (cf. Pruss \cite{prus})
uses in the same context is much stronger and obviously  implies weak mean bounded condition.

In Section 4 we will consider martingale random field, thus introduce the fundamental notions.
Let $\{\mathfrak{F}_{{\bf k}}, \enskip  {\bf k} \in \mathbb{N}^d \}$ be a filtration of $\sigma$-algebras i.e.
\begin{description}
  \item[(F1)] \enskip if ${\bf k}\leq {\bf n}$ \enskip $\Rightarrow$ \enskip $ {\mathfrak{F}}_{\bf k}\subset {\mathfrak{F}}_{\bf n}\subset {\mathfrak{F}}$
\end{description}

\medskip
An integrable family of random variables
$ \{ Z_{\bf k}, {\bf k}\in \mathbb{N}^d\}$,  adapted to $\{\mathfrak{F}_{\bf k}, {\bf k} \in \mathbb{N}^d \}$   is called
martingale random field,  iff
$$\bigwedge \limits_{{\bf k}\leq {\bf n}} \enskip
 E\bigl(Z_{\bf n} \mid  \mathfrak{F}_{\bf k}\bigr) = Z_{{\bf k}} \quad \text{ a.s. }
$$
Let us observe, that for martingale $\{(Z_{\bf n}, \mathfrak{F}_{\bf n}), {\bf n}\in \mathbb{N}^d\}$
$$
X_{\bf n} = \sum_{{\bf a}\in \{0, 1\}^r} (-1)^{\sum_{i=1}^r a_i} Z_{{\bf n}- {\bf a}},
$$
where ${\bf a} = (a_1, a_2, \dots , a_r)$ and ${\bf n}\in \mathbb{N}^d$,
are martingale differences with respect to $\{\mathfrak{F}_{\bf n}, {\bf n}\in \mathbb{N}^d\}$.

\bigskip
\section{ Auxiliary Lemmas }

\bigskip
Let $\{ a_{\mathbf{k},\mathbf{n}},
\mathbf{k},\mathbf{n} \in  \mathbb{N}^d \}$ be a family of real numbers , such that
$ 0 \leq a_{\mathbf{k},\mathbf{n}} < 1 $, then we have.
\bigskip

\noindent
\begin{lemma} \label{expineq} If $\sum \limits _{\mathbf{k} \leq \mathbf{n}} a_{\mathbf{k},\mathbf{n}}\rightarrow 0 $ as
 $\mathbf{n} \rightarrow \infty$, then for a given $0 < \delta < 1$ and $\mathbf{n}$ sufficiently large
$$ 1-\prod  \limits _{\mathbf{k} \leq \mathbf{n}}(1- a_{\mathbf{k},\mathbf{n}}) \geq
(1-\delta)\sum \limits _{\mathbf{k} \leq \mathbf{n}} a_{\mathbf{k},\mathbf{n}}.$$
\end{lemma}

\begin{pr} Let, for a given $\delta$,  $\mathbf{n}$ be sufficiently large such that
$\sum \limits _{\mathbf{k} \leq \mathbf{n}} a_{\mathbf{k},\mathbf{n}} \leq \delta(1- \delta)$.\\
Thus
\begin{equation}
\begin {split}
&\prod  \limits _{\mathbf{k} \leq \mathbf{n}}(1- a_{\mathbf{k},\mathbf{n}})  =
  \exp\{\sum \limits _{\mathbf{k} \leq
   \mathbf{n}}\ln(1- a_{\mathbf{k},\mathbf{n}})\} \leq
    \exp\{-\sum \limits _{\mathbf{k} \leq
     \mathbf{n}} a_{\mathbf{k},\mathbf{n}}\} \leq \\&
     1- \sum \limits _{\mathbf{k}  \mathbf{n}} a_{\mathbf{k},\mathbf{n}} +( \sum \limits _{\mathbf{k} \leq
     \mathbf{n}} a_{\mathbf{k},\mathbf{n}})^2 \leq
     1-(1-\delta)\sum \limits _{\mathbf{k} \leq
     \mathbf{n}} a_{\mathbf{k},\mathbf{n}}
     \end{split}
     \end{equation}
Now, assertion easily follows.
\end{pr}

The next lemma is simply consequence of WMB  condition and the well known fact, that
for any random variable $X$ with $E\mid X\mid ^s < \infty$
$$
E\mid X\mid ^s = s \int_0^\infty x^{s-1}P[\mid X\mid >x]dx.
$$
For some $a>0,$ let us put
$$
X_{\bf i}' = X_{\bf i}I[\mid X_{\bf i}\mid \leq a],  \quad X_{\bf i}'' =X_{\bf i} I[\mid X_{\bf i}\mid >a ],
$$
and
$$
\xi' = \xi I[\mid \xi\mid \leq a],  \quad \xi'' =\xi I[\mid \xi\mid >a ].
$$

\medskip
\begin{lemma}  \label{nierowmomen}
Let $\{X_{\bf n}, {\bf n}\in \mathbb{N}^d\}$ be a field of random variables satisfying
WMB  condition with  random variable $\xi$ and constants $\kappa_1,\kappa_2$.
Let $s>0$.

\medskip
(a) If $E\mid \xi\mid ^s <\infty$, then $\kappa_2 E\mid \xi\mid ^s \leq {1\over \mid {\bf n}\mid }\sum_{{\bf k}\leq {\bf n}}E\mid X_{\bf k}\mid ^s
\leq \kappa_1 E\mid \xi\mid ^s$.

\medskip
(b) $\kappa_2 E\mid \xi'\mid ^s \leq {1\over \mid {\bf n}\mid }\sum_{{\bf k}\leq {\bf n}}E\mid X_{\bf k}'\mid ^s \leq \kappa_1 E\mid \xi'\mid ^s  $ .

\medskip
(c) $\kappa_2E\mid \xi''\mid ^s \leq {1\over \mid {\bf n}\mid }\sum_{{\bf k}\leq {\bf n}}E\mid X_{\bf k}''\mid ^s \leq \kappa_1E\mid \xi''\mid ^s $.

\end{lemma}
\medskip
The following properties of ND random variables, proved by Bozorgnia et al. \cite{bozor}, for sequences of r.v., obviously
hold true for ND random fields.
\begin{lemma} \label{w³asnoND} Let
$\{ X_{\bf k}, {\bf k} \leq {\bf n} \}$ be a field of  ND random
variables and $\{ f_{\bf k}, {\bf k} \leq {\bf n} \}$ a family of Borel functions,  which all are non-decreasing (non-increasing),
then\\
(a)  $\{ f(X_{\bf k}), {\bf k} \leq {\bf n} \}$ is a ND random field,\\
\medskip
(b)  if  additionally, $X_{\bf k}$ are non-negative, we have $$E\left(\prod \limits_{{\bf k} \leq {\bf n}}X_{\bf k}\right) \leq \prod \limits_{{\bf k} \leq {\bf n}}
EX_{\bf k}.$$

\end{lemma}
\medskip
\begin{lemma}\label{expobciete}
Assume, that $\{ X_{\bf n}, {\bf n}\in \mathbb{N}^d \}$ is a field of zero mean, square integrable   ND random
variables WMD by random variable $\xi$ and such that $E\xi^2 = \sigma^2< \infty $, then\\
(a) \hspace{18mm}  $E(\sum \limits_ {{\bf k}\leq {\bf n}}X_{\bf k})^2  \leq \kappa_1 \sigma^2  \mid\mathbf{n}\mid,$\\
  if additionally  $P(X_{\bf k}\leq b)=1$ for every  ${\bf k}\leq {\bf n},$ then \\
(b) \hspace{18mm} $P(\sum \limits_ {{\bf k}\leq {\bf n}}X_{\bf k}>x) \leq
e^{-tx +\kappa_1 \sigma^2  \mid\mathbf{n}\mid}$  \\
 for all $x,b>0$ and $0<t<\frac{1}{b}.$

\end{lemma}
\begin{pr}
By lemma \ref{w³asnoND} and  \ref{nierowmomen}a one can obtain
$$E(\sum \limits_ {{\bf k}\leq {\bf n}}X_{\bf k})^2 \leq  \sum \limits_ {{\bf k}\leq {\bf n}}EX_{\bf k}^2
\leq \kappa_1 \sigma^2  \mid\mathbf{n}\mid,$$
thus (a) holds. Standard  inequalities: $e^x \leq 1+x+x^2$ for $0<x<1$, $1+x<e^x$ Lemma \ref{w³asnoND}
and (a) lead us to inequality (b).
\end{pr}

\medskip
Let us put $M_{\mathbf{n}}^r=\sum \limits _{\mathbf{k} \leq \mathbf{n}}E\mid X_{\mathbf{k}}\mid ^r$, $ \lambda_{\mathbf{k}}=
E X_{\mathbf{k}}I[X_{\mathbf{k}} \geq -y]$ and $\Lambda_{\mathbf{n}}=\sum \limits _{\mathbf{k} \leq \mathbf{n}}
\lambda_{\mathbf{k},y}$

\begin{lemma}\label{hofjond} Assume that $\{ X_{\bf n}, {\bf n}\in \mathbb{N}^d \}$ be a field of zero mean ND random
variables with finite an absolute r-th moment, $1 \leq r \leq 2$, then exist constant $C>0$ such that for every $x>0$ and $j > 0$
$$P(\mid S_{\mathbf{n}}\mid >x) \leq P(\max \limits_{\mathbf{k} \leq \mathbf{n}}\mid X_{\mathbf{k}}\mid >\frac{x}{j})+C(\frac{1}{x^r}
M_{\mathbf{n}}^r)^j$$
\end{lemma}

\begin{pr} Fakoor et al. have proved Fuk-Nagaev inequality for sequences ND random variables,  Theorem 3 of \cite{Fako},
since the proof  doesn't involve the order of index set, inequality holds true for $d \geq 2$ case, thus under assumption  for any $y>0$ we have
\begin{equation}
\begin{split}
&P(\mid S_{\mathbf{n}}\mid >x) \leq P(\max \limits_{\mathbf{j} \leq \mathbf{n}}\mid X_{\mathbf{j}}\mid >y)+
\\& 2\exp\{\frac {x}{y}-\left(\frac {x}{y}-\frac{\Lambda_{\mathbf{n}}}{y}
    +\frac{M_{\mathbf{n}}^r}{y^r}\right)\ln (1+\frac{xy^{r-1}}{M_{\mathbf{n}}^r})\}=I_1
\end{split}
\end{equation}

Since for all ${\bf k}\in \mathbb{N}^d $
$$\frac{\lambda_{\mathbf{k},y}}{y}=\frac{E \mid X_{\mathbf{k}}\mid I[\mid X_{\mathbf{k}}\mid  \geq y]}{y} \leq
 \frac{E \mid X_{\mathbf{k}}\mid ^r}{y^r},$$

thus putting $\frac{x}{y}=j, $ we obtain

\begin{equation*}
\begin{split}
I_1 \leq &P(\max \limits_{\mathbf{k} \leq \mathbf{n}}\mid X_{\mathbf{k}}\mid >\frac{x}{j})+
2\exp\Big\{j-j\ln\left(1+\frac{x^rj^{1-r}}
{M_{\mathbf{n}}^r}\right)\Big\}
\leq  \\&P(\max \limits_{\mathbf{k} \leq \mathbf{n}}\mid X_{\mathbf{k}}\mid >
\frac{x}{j})+ 2e^jj^{(r-1)j}\left(\frac{M_{\mathbf{n}}^r}
{x^r}\right)^j,
\end{split}
\end{equation*}
which finishes the proof of lemma.
\end{pr}

\begin{lemma}\label{sumkwad}
 Let $\xi$ be a random variable such that $E\mid \xi\mid ^{\frac{1}{\alpha_1}}(\log_{+}\mid \xi\mid )^{p-1} < \infty$,
then under our setting with $\alpha_1>\frac{1}{2}$
\begin{equation*}
\sum\limits _{\mathbf{n} \in \mathbb{N}^d}\frac{E\mid \xi\mid ^2I[\mid \xi\mid  \leq
\mid\mathbf{n}^{\boldsymbol{\alpha}}\mid]}{\mid\mathbf{n}^{\boldsymbol{\alpha}}\mid^2} < \infty.
\end{equation*}
\end{lemma}
\begin{pr} For every $\nu \in \mathbf{N}$ define
$$\Delta f(\nu)=\text{card} \{(n_1,n_2,\ldots, n_p): n_1\cdot n_2\cdot\ldots\cdot n_p= \nu\}.$$
From the proof of Theorem 2.1 by Gut et al. \cite{Gut Asy} one can deduce
\begin{equation*}
\begin{split}
&\sum\limits _{\mathbf{n} \in \mathbb{N}^d}\frac{E\mid \xi\mid ^2I[\mid \xi\mid  \leq
\mid\mathbf{n}^{\boldsymbol{\alpha}}\mid]}{\mid\mathbf{n}^{\boldsymbol{\alpha}}\mid^2} =
\\&
\sum\limits_{\nu=1}^{\infty}
\sum \limits_{n_{p+1},\ldots,n_d=1}^{\infty}\Delta f(\nu)\frac{1}{\nu^{2\alpha_1}\cdot n_{p+1}^{2\alpha_{p+1}}\cdot \ldots
n_d^{2\alpha_d}}\sum\limits_{j=1}^{\nu^{\alpha_1}\cdot n_{p+1}^{\alpha_{p+1}}\cdot \ldots \cdot
n_d^{\alpha_d}}E\xi^2I[j-1< \xi \leq j] \leq
\\&CE\mid \xi\mid ^{\frac{1}{\alpha_1}}(\log_{+}\mid \xi\mid )^{p-1} < \infty,
\end{split}
\end{equation*}
where $C>0$ is suitable constant.
\end{pr}

\section{Baum-Katz  type theorems for ND  random fields}

The first two theorems of this Section are extensions  and  compliments of some results of  Peligrad \cite{pelgut},
Gut et al.\cite{GutHsu}, \cite{Gut Asy},  Kuczmaszewska et al. \cite {Kucz}.
\begin{theorem}
Let  $r \geq 1$, $\alpha_1\geq \frac{1}{2}$,  $\alpha_1r \geq 1$ and $\{ X_{\bf n}, {\bf n}\in \mathbb{N}^d \}$ be a
zero mean  random field of
 ND random variables, weak mean bounded by $\xi$. If
\begin{equation} \label{moment}
E\mid \xi\mid ^r (\log_{+}\mid \xi\mid) ^{p-1}<\infty,
\end{equation}
then
\begin{equation} \label{hsu}
\sum\limits_{\mathbf{n}}\mid \mathbf{n}\mid ^{\alpha_1r-2}P(\mid S_{\mathbf{n}}\mid >\mid \mathbf{n}^{\boldsymbol{\alpha}}\mid \varepsilon)< \infty   \quad \text{for all}\quad  \varepsilon >0.
\end{equation}
Conversly if
\begin{equation}\label{hsumax}
\sum\limits_{\mathbf{n}}\mid \mathbf{n}\mid ^{\alpha_1r-2}
P(\max\limits_{\mathbf{k}\leq\mathbf{n}}\mid S_{\mathbf{k}}\mid >\mid \mathbf{n}^{\boldsymbol{\alpha}}\mid \varepsilon)<
\infty  \quad \text{for all}\quad \varepsilon >0,
\end{equation}
then (\ref{moment}) holds.

\end{theorem}

\begin{pr} $(\ref{moment}) \Rightarrow (\ref{hsu})$.
The general idea of the proof is based on the proof of Theorem 4.1 by Gut and Stadm\"{u}ller \cite{GutHsu},
 thus we sketch the proof showing differences. At the beginning, assume that  $\alpha_1> \frac{1}{2}$,  $\alpha_1r > 1$ and (\ref{moment}) holds.
 Applying Lemma \ref{hofjond} one can obtain
\begin{equation}  \begin{split}
&\sum\limits_{\mathbf{n}}\mid  \mathbf{n}\mid  ^{\alpha_1r-2}P(\mid  S_{\mathbf{n}}\mid  >\mid  \mathbf{n}^{\boldsymbol{\alpha}}\mid   \varepsilon) \leq
\sum \limits_{\mathbf{n}} \mid  \mathbf{n}\mid  ^{\alpha_1r-2} \sum \limits_{\mathbf{k} \leq \mathbf{n}} P(\mid  X_{\mathbf{k}}\mid  >y)+\\&
+\frac{C}{\varepsilon^{rj}}\sum\limits_{\mathbf{n}}\mid  \mathbf{n}\mid  ^{\alpha_1r-2}\mid  \mathbf{n}^{\boldsymbol{\alpha}}\mid   ^{-jr} (\sum \limits_{\mathbf{k} \leq \mathbf{n}}
E\mid  X_{\mathbf{k}}\mid  ^r)^j \leq \\&
C_1\sum\limits_{\mathbf{n}}\mid  \mathbf{n}\mid  ^{\alpha_1r-1}
P(\mid  \xi\mid  >\mid  \mathbf{n}^{\boldsymbol{\alpha}}\mid  \varepsilon^{'})+ C_2\sum\limits_{\mathbf{n}}\mid  \mathbf{n}\mid  ^{\alpha_1r-2+j}\mid  \mathbf{n}^{\boldsymbol{\alpha}}\mid   ^{-jr}
(E\mid  \xi\mid  ^r)^j\\&
=I_2+I_3, \quad \text{where}\quad \varepsilon^{'}=\frac{\varepsilon}{j}.
\end{split}
\end{equation}
\noindent
The first sum  $I_2$  is finite by Lemma 2.2 of \cite{GutHsu}, the second one is estimated as follows

\begin{equation} \begin{split}
&I_2\leq C\sum\limits_{\mathbf{n}}\mid  \mathbf{n}\mid  ^{\alpha_1r-2+j}\mid  \mathbf{n}^{\boldsymbol{\alpha}}\mid   ^{-jr} \leq \\&
C\sum\limits_{\mathbf{n}}\prod\limits_{i=1}^dn_i^{\alpha_1r-2+(1-\alpha_1r)j}\leq
C\prod\limits_{i=1}^d \sum\limits_{n_i=1}^{\infty} n_i^{\alpha_1r-2+(1-\alpha_1r)j}<\infty,
\end{split}
\end{equation}

since exponent in the last sum can be less than $(-1)$, for j sufficiently large.

Now, assume that $\alpha_1>\frac{1}{2}$, $\alpha_1r=1.$
Let $Y_{\mathbf{k},\mathbf{n}}=\min(\mid  \mathbf{n}\mid  ^{\boldsymbol{\alpha}},\mid X_{\mathbf{k}}\mid )\sgn (X_{\mathbf{k}})$,
$X_{\mathbf{k},\mathbf{n}}=X_{\mathbf{k}}I[\mid X_{\mathbf{k}}\mid   \leq \mid  \mathbf{n}\mid  ^{\boldsymbol{\alpha}}]$
and $T_{\mathbf{n}}=\sum \limits_{\mathbf{k} \leq \mathbf{n}}Y_{\mathbf{k},\mathbf{n}}.$ Thus we get

\begin{equation} \begin{split}
&\sum\limits_{\mathbf{n}}\frac{1}{\mid  \mathbf{n}\mid  }P(\mid  S_{\mathbf{n}}\mid  >2\mid  \mathbf{n}^{\boldsymbol{\alpha}}\mid  \varepsilon) \leq \\&
\sum\limits_{\mathbf{n}}\frac{1}{\mid  \mathbf{n}\mid  }P(\mid  T_{\mathbf{n}}\mid  >\mid  \mathbf{n}^{\boldsymbol{\alpha}}\mid  \varepsilon)+
\sum\limits_{\mathbf{n}}\frac{1}{\mid  \mathbf{n}\mid  }P(\mid  S_{\mathbf{n}}-T_{\mathbf{n}}\mid  >\mid  \mathbf{n}^{\boldsymbol{\alpha}}\mid
\varepsilon)=I_4+I_5
\end{split}
\end{equation}
The first sum can be estimated by applying Chebyshev inequality, Lemma \ref{expobciete} and \ref{nierowmomen},  WMD condition consecutively:
\begin{equation}\begin{split}
&I_4\leq\sum\limits_{\mathbf{n}}\frac{1}{\mid  \mathbf{n}\mid  }\frac{E(T_{\mathbf{n}}-ET_{\mathbf{n}})^2}{\varepsilon^2\mid  \mathbf{n}^{\boldsymbol{\alpha}}\mid  ^2}
\leq
C\sum\limits_{\mathbf{n}}\frac{1}{\mid  \mathbf{n}\mid  }\frac{ET_{\mathbf{n}}^2}{\mid  \mathbf{n}^{\boldsymbol{\alpha}}\mid  ^2}
\leq \\&
C(\sum\limits_{\mathbf{n}}[ \frac{1}{\mid  \mathbf{n}\mid  }
\frac{\sum\limits_{\mathbf{k}\leq \mathbf{n}}EX_{\mathbf{k},\mathbf{n}}^2}{\mid  \mathbf{n}^{\boldsymbol{\alpha}}\mid  ^2}+
\frac{1}{\mid  \mathbf{n}\mid  }\sum\limits_{\mathbf{k}\leq \mathbf{n}}P(\mid  X_{\mathbf{k}}\mid  >\mid  \mathbf{n}^{\boldsymbol{\alpha}})])\leq \\&
C(\sum\limits_{\mathbf{n}}\frac{E\mid  \xi\mid  ^2I[\mid  \xi\mid   \leq
\mid \mathbf{n}^{\boldsymbol{\alpha}}\mid ]}{\mid \mathbf{n}^{\boldsymbol{\alpha}}\mid^2}+
\sum\limits_{ \mathbf{n}}P(\mid  \xi\mid  >\mathbf{n}^{\boldsymbol{\alpha}}))\leq CE\mid  \xi\mid  ^{\frac{1}{\alpha_1}}(\log_{+}\mid  \xi\mid  )^{p-1}.
\end{split}
\end{equation}
The last inequality follows from Lemma \ref{sumkwad} and Lemma 2.2 of \cite{GutHsu} respectively.\\
On the other hand

\begin{equation}\begin{split}
&I_5\leq\sum\limits_{\mathbf{n}}\frac{1}{\mid \mathbf{n}\mid}
P(\sum\limits_{\mathbf{k}\leq \mathbf{n}}
\mid X_{\mathbf{k}}\mid I[\mid  X_{\mathbf{k}}\mid >\mid \mathbf{n}^{\boldsymbol{\alpha}}\mid ]>
\varepsilon \mid \mathbf{n}^{\boldsymbol{\alpha}}\mid )
\leq\\&
\sum\limits_{\mathbf{n}}\frac{1}{\mid \mathbf{n}\mid }
P(\sum\limits_{\mathbf{k}\leq \mathbf{n}}
\mid X_{\mathbf{k}}\mid > \mid \mathbf{n}^{\boldsymbol{\alpha}}\mid )\leq C\sum\limits_{\mathbf{n}}
P(
\mid \xi\mid > \mid \mathbf{n}^{\boldsymbol{\alpha}}\mid )<\infty,
\end{split}
\end{equation}
by WMD condition and Lemma 2.2 of \cite{GutHsu}.\\
The implication (\ref{hsumax})$\Rightarrow$ (\ref{moment}). Firstly, let us observe, that the negative and positive part
of ND random variables are still ND. Thus

\begin{equation}\begin{split}
&P(\max\limits_{\mathbf{k}\leq\mathbf{n}}\mid S_{\mathbf{k}}\mid>\mid \mathbf{n}^{\boldsymbol{\alpha}}\mid \varepsilon)\geq
P(\max\limits_{\mathbf{k}\leq\mathbf{n}}\mid X_{\mathbf{k}}\mid>2\mid \mathbf{n}^{\boldsymbol{\alpha}}\mid \varepsilon)\geq \\&
P(\max\limits_{\mathbf{k}\leq\mathbf{n}}X^{+}_{\mathbf{k}}>2\mid \mathbf{n}^{\boldsymbol{\alpha}}\mid \varepsilon)=
1-P(\bigcap\limits_{\mathbf{k}\leq\mathbf{n}}[X^{+}_{\mathbf{k}}\leq2\mid \mathbf{n}^{\boldsymbol{\alpha}}\mid \varepsilon])\geq \\&
1-\prod \limits_{\mathbf{k}\leq\mathbf{n}}P(X^{+}_{\mathbf{k}}\leq2\mid \mathbf{n}^{\boldsymbol{\alpha}}\mid \varepsilon) =
1-\prod \limits_{\mathbf{k}\leq\mathbf{n}}(1-P(X^{+}_{\mathbf{k}}>2\mid \mathbf{n}^{\boldsymbol{\alpha}}\mid \varepsilon))
\end{split}
\end{equation}
From  (\ref{hsumax})  and  (15) it's easy to see, that $\prod \limits_{\mathbf{k}\leq\mathbf{n}}(1-
P(X^{+}_{\mathbf{k}}>2\varepsilon\mid \mathbf{n}^{\boldsymbol{\alpha}}\mid )) \rightarrow 1$ as $\mathbf{n} \rightarrow \infty,$
what is equivalent to
\begin{equation}
\sum\limits_{\mathbf{k}\leq\mathbf{n}}
P(X^{+}_{\mathbf{k}}>2\varepsilon\mid \mathbf{n}^{\boldsymbol{\alpha}}\mid )\rightarrow 0
\quad \text{as}\quad  \mathbf{n} \rightarrow \infty.
\end{equation}
Analogously, we can get

\begin{equation}
\sum\limits_{\mathbf{k}\leq\mathbf{n}}
P(X^{-}_{\mathbf{k}}>2\varepsilon\mid \mathbf{n}^{\boldsymbol{\alpha}}\mid )\rightarrow 0
\quad \text{as}\quad  \mathbf{n} \rightarrow \infty.
\end{equation}
Now, applying Lemma \ref{expineq} with $a_{\mathbf{k},\mathbf{n}}=\sum\limits_{\mathbf{k}\leq\mathbf{n}}
P(  X_{\mathbf{k}}^{+}  >\varepsilon\mid \mathbf{n}^{\boldsymbol{\alpha}}\mid )$ and \\$a_{\mathbf{k},\mathbf{n}}=\sum\limits_{\mathbf{k}\leq\mathbf{n}}
P(  X_{\mathbf{k}}^{-}  >\varepsilon\mid \mathbf{n}^{\boldsymbol{\alpha}}\mid ),$
 WMB condition and Lemma 2.2 of \cite{GutHsu},\\ we have

\begin{equation}\begin{split}
&\sum\limits_{\mathbf{n}}\mid \mathbf{n}\mid ^{\alpha_1r-2}
P(\max\limits_{\mathbf{k}\leq\mathbf{n}}\mid S_{\mathbf{k}}\mid >\varepsilon \mid \mathbf{n}^{\boldsymbol{\alpha}}\mid ) \geq
C_1\sum\limits_{\mathbf{n}}\mid \mathbf{n}\mid ^{\alpha_1r-2}
P(\max\limits_{\mathbf{k}\leq\mathbf{n}}\mid X_{\mathbf{k}}\mid >2\varepsilon\mid  \mathbf{n}^{\boldsymbol{\alpha}}\mid )\geq \\&
C_2\sum\limits_{\mathbf{n}}\mid \mathbf{n}\mid ^{\alpha_1r-2}
\sum \limits_{\mathbf{k}\leq\mathbf{n}}P(\mid X_{\mathbf{k}}\mid >2\varepsilon\mid  \mathbf{n}^{\boldsymbol{\alpha}}\mid )\geq
C_3\sum\limits_{\mathbf{n}}\mid \mathbf{n}\mid ^{\alpha_1r-1}
\sum \limits_{\mathbf{k}\leq\mathbf{n}}P(\mid \xi\mid >2\varepsilon\mid  \mathbf{n}^{\boldsymbol{\alpha}}\mid )\geq \\&
C_4E\mid \xi\mid ^r (\log_{+}\mid \xi\mid)\mid ^{p-1}.
\end{split}
\end{equation}
\end{pr}

\begin{theorem} Let  $\{ X_{\bf n}, {\bf n}\in \mathbb{N}^d \}$ be a field of zero mean ND random variables
satisfying WMB condition with r.v. $\xi$ and suppose, that $r\geq 2$, $\alpha_1= \frac{1}{2}$,  $\alpha_1r \geq 1.$
If
\begin{equation}
\label{momentlog}
E\mid \xi\mid ^r (\log_{+}\mid \xi\mid)\mid ^{p-1-\frac{r}{2}}<\infty \quad \text{and} \quad E\xi^2=\sigma^2 < \infty,
\end{equation}
then
\begin{equation}\label{hsulog}
\sum\limits_{\mathbf{n}}\mid \mathbf{n}\mid ^{(r/2)-2}P(\mid S_{\mathbf{n}}\mid \geq
\sqrt{ \prod\limits_{i=1}^pn_i\log(\prod\limits_{i=1}^pn_i)}\prod\limits_{i=p+1}^dn_i^{\alpha_i}\varepsilon)< \infty
\end{equation}
for $ \varepsilon>\sigma_1 \sqrt{r-2},$
where $p=\max\{k:\alpha_k=\alpha_1\}$  and $\sigma_1^2=\kappa_1\sigma^2$.\\
Conversely, suppose either $r=2$ and $p \geq 2$ or that  $r> 2$.\\
If
\begin{equation}\label{hsulogmax}
\sum\limits_{\mathbf{n}}\mid \mathbf{n}\mid ^{(r/2)-2}P(\max \limits_{\mathbf{k} \leq \mathbf{n}}\mid S_{\mathbf{n}}\mid \geq
\sqrt{ \prod\limits_{i=1}^pn_i\log(\prod\limits_{i=1}^pn_i)}\prod\limits_{i=p+1}^dn_i^{\alpha_i}\varepsilon) < \infty
\end{equation}
for some $ \varepsilon >0,$
then
\begin{equation} \label{momentlog2}
E\mid \xi\mid ^r (\log_{+}\mid \xi\mid)\mid ^{p-1-\frac{r}{2}}<\infty.
\end{equation}
\end{theorem}
 \begin{pr} (\ref{momentlog} ) $\Rightarrow$ (\ref{hsulog}), the case $\alpha_1=\frac{1}{2},$ $r=2.$\\
 Applying Lemma \ref{hofjond} with second moment and by WMB condition, we have

 \begin{equation}
 \begin{split}
& \sum\limits_{\mathbf{n}}\frac{1}{\mid \mathbf{n}\mid }P\left( \mid S_{\mathbf{n}}\mid >
\sqrt{ \prod\limits_{i=1}^pn_i\log(\prod\limits_{i=1}^pn_i)}\prod\limits_{i=p+1}^dn_i^{\alpha_i}\varepsilon \right) \leq\\&
2\sum\limits_{\mathbf{n}}\frac{1}{\mid \mathbf{n}\mid}\sum\limits_{\mathbf{k}\leq\mathbf{n}}P\left(\mid X_{\mathbf{k}}\mid>
\sqrt{\prod\limits_{i=1}^pn_i\log(\prod\limits_{i=1}^pn_i)}\prod\limits_{i=p+1}^dn_i^{\alpha_i}\varepsilon^{'}\right)+\\&
\frac{C}{\varepsilon^2j} \sum\limits_{\mathbf{n}}\frac{1}{\mid \mathbf{n}\mid}
\left(\frac{\sum\limits_{\mathbf{k}\leq\mathbf{n}}EX^2_{\mathbf{k}}}
{ \prod\limits_{i=1}^pn_i\log(\prod\limits_{i=1}^pn_i)(\prod\limits_{i=p+1}^dn_i^{\alpha_i})^2}\right)^j \leq
\end{split}
 \end{equation}
 \begin{equation*}
 \begin{split}
&C_1\sum\limits_{\mathbf{n}}P\left(\mid \xi \mid >
\sqrt{ \prod\limits_{i=1}^pn_i\log(\prod\limits_{i=1}^pn_i)}\prod\limits_{i=p+1}^dn_i^{\alpha_i}\varepsilon^{'}\right)+\\
&C_2\sum\limits_{\mathbf{n}}\frac{1}{\mid\mathbf{n}\mid}
\left(\frac{\mid\mathbf{n}\mid  E\xi^2}
{\prod\limits_{i=1}^pn_i\log(\prod\limits_{i=1}^pn_i)\prod\limits_{i=p+1}^dn_i^{2\alpha_i}}\right)^j=I_4+I_5,
\end{split}
 \end{equation*}
 where $C_1$ and $C_2$ are suitable constants. The first sum is finite by assumption (\ref{momentlog})
 and the second one, by the same arguments as in proof of Theorem 4.1 of \cite{GutHsu}.  \\
 \medskip
The case $\alpha_1=\frac{1}{2}, \enskip r>2.$ Let $0< \eta < \alpha_{p+1}- \frac{1}{2}$ and $\beta _i =\alpha_i - \eta$
for \\
$i=p+1,p+2, \ldots ,d$. Furthermore, we set
\begin{equation*}
\begin{split}
&a_{\mathbf{n}}=\sqrt{ \prod\limits_{i=1}^pn_i\log(\prod\limits_{i=1}^pn_i)}\prod\limits_{i=p+1}^dn_i^{\alpha_i},  \enskip
b_{\mathbf{n}}=\frac{2\delta\sigma^2_1}{\varepsilon} \sqrt{ \frac{\prod\limits_{i=1}^pn_i}{\log(\prod\limits_{i=1}^pn_i)}}
\prod\limits_{i=p+1}^dn_i^{\beta_i}   \\&
c_{\mathbf{n}}=\delta a_{\mathbf{n}} \enskip \text{and} \enskip
d_{\mathbf{n}}=\sqrt{ \prod\limits_{i=1}^pn_i \log(\prod\limits_{i=1}^pn_i)}\prod\limits_{i=p+1}^dn_i^{1-\beta_i}
 \end{split}
  \end{equation*}

\noindent
 Let us put $X_{\mathbf{k}}^{'}=X_{\mathbf{k}}I[X_{\mathbf{k}} < b_{\mathbf{n}}] +
 b_{\mathbf{n}}I[X_{\mathbf{k}} \geq b_{\mathbf{n}}]$ and  $S_{\mathbf{k}}^{'}=\sum\limits_{\mathbf{k}\leq\mathbf{n}}
 X_{\mathbf{k}}^{'}.$\\
   Define the events:
\begin{equation*}
\begin{split}
&A^1_{\mathbf{n}}=\{ S_{\mathbf{n}}^{'}>\varepsilon a_{\mathbf{n}}\}, \enskip
  A^2_{\mathbf{n}}=\{ \text{at least two} \enskip   \mathbf{k}, \enskip \mathbf{k}\leq\mathbf{n}:  b_{\mathbf{n}}<
  X_{\mathbf{k}}<c_{\mathbf{n}}\}\\&
  A^3_{\mathbf{n}}=\{ \text{at least one} \enskip   \mathbf{k}, \enskip \mathbf{k}\leq\mathbf{n}:
  X_{\mathbf{k}} \geq c_{\mathbf{n}}\}. \enskip \text{and}  \enskip A_{\mathbf{n}}=\{ S_{\mathbf{n}}>
  (\varepsilon +2\delta) a_{\mathbf{n}}\}.
  \end{split}
  \end{equation*}
  It's clearly that $A_{\mathbf{n}} \subset A^1_{\mathbf{n}} \cup A^2_{\mathbf{n}} \cup A^3_{\mathbf{n}}$ thus
  $P(A_{\mathbf{n}}) \leq P( A^1_{\mathbf{n}}) +P(A^2_{\mathbf{n}})+P( A^3_{\mathbf{n}}).$\\
  We start from the estimation of $P( A^1_{\mathbf{n}}).$
  The first step, since $\{ X_{\mathbf{ k}}, \mathbf{ k}\leq  \mathbf{ n} \}$
  is a field of zero mean random variables satisfying WMB condition, we have
 \begin{equation*}
\begin{split}
& \mid ES^{'}_{\mathbf{n}}\mid \leq
\sum  \limits _{\mathbf{ k}\leq  \mathbf{ n}}E X_{\mathbf{ k}}
I[X_{\mathbf{k}} \geq b_{\mathbf{n}}] \leq \kappa_1\mid \mathbf{n} \mid E\xi I[\xi \geq b_{\mathbf{n}}]\leq \\&
\frac{\kappa_1\mid \mathbf{n} \mid}{b_{\mathbf{n}}}E\xi^2 I[\xi \geq b_{\mathbf{n}}]=o(d_{\mathbf{n}})
 \end{split}
  \end{equation*}
  Further arguments and details are the same as proof of (4.4-4.6) of \cite{GutHsu}, hence
   \begin{equation}\label{prA1}
  \sum  \limits _{ \mathbf{ n}} \mid  \mathbf{ n} \mid ^{(r/2-2)} P( A^1_{\mathbf{n}})< \infty \quad \text{for}\quad
  \varepsilon > \sigma_{1}\frac{1+\delta}{1-\delta}\sqrt{r-2} \quad \text{and all }\quad \delta>0.
    \end{equation}
  In estimation of $P( A^2_{\mathbf{n}}),$  we exploit the ND and WMD property of
   $\{ X_{\mathbf{ k}}, \mathbf{ k}\leq  \mathbf{ n} \}$ and thereafter by the same manner as in the proof of (4.7)
   of \cite{GutHsu}
    \begin{equation*}
    P( A^2_{\mathbf{n}}) \leq\sum\limits_{\mathbf{ k}\leq  \mathbf{ n}} \sum \limits_{\mathbf{ l}\leq  \mathbf{ n},
    \mathbf{ l}\neq \mathbf{ k}} P(X_{\mathbf{ k}}>b_{\mathbf{n}},X_{\mathbf{ l}}>b_{\mathbf{n}}) \leq
    \kappa^2_1 \mid  \mathbf{ n} \mid^2(P( \xi > b_{\mathbf{n}}))^2
     \end{equation*}
     thus

     \begin{equation} \label{prA2}
      \begin{split}
    &  \sum  \limits _{ \mathbf{ n}} \mid  \mathbf{ n} \mid ^{(r/2-2)}  P( A^2_{\mathbf{n}}) \leq \\&
      C(r,\delta)
       \sum  \limits _{ \mathbf{ n}} \mid  \mathbf{ n} \mid ^{-r/2}
       \frac{\left(\log(\prod\limits_{i=1}^pn_i)\right)^r}{\left(\log \mid  \mathbf{ n} \mid\right)^{2(p-1)-r}}< \infty
        \quad \text{for all }\quad \delta>0.
     \end{split}
     \end{equation}
    Finally, by Lemma 2.1(c) of \cite{GutHsu}
       \begin{equation} \label{prA3}
       \begin{split}
      & \sum  \limits _{ \mathbf{ n}} \mid  \mathbf{ n} \mid ^{(r/2-2)}  P( A^3_{\mathbf{n}}) \leq
       \sum  \limits _{ \mathbf{ n}} \mid  \mathbf{ n} \mid ^{(r/2-2)} \sum\limits_{\mathbf{ k}\leq  \mathbf{ n}}
       P( X_{\mathbf{k}} \geq c_{\mathbf{n}}) \leq \\&
        \kappa_1 \sum  \limits _{ \mathbf{ n}} \mid  \mathbf{ n} \mid ^{(r/2-1)}
        P( \xi \geq \delta
    \sqrt{ \prod\limits_{i=1}^pn_i\log(\prod\limits_{i=1}^pn_i)}\prod\limits_{i=p+1}^dn_i^{\alpha_i} ) < \infty
     \end{split}     \end{equation}

Now,  let us put $X_{\mathbf{k}}^{''}=X_{\mathbf{k}}I[X_{\mathbf{k}}>- b_{\mathbf{n}}] -
 b_{\mathbf{n}}I[X_{\mathbf{k}} \leq -b_{\mathbf{n}}]$,  $S_{\mathbf{k}}^{''}=\sum\limits_{\mathbf{k}\leq\mathbf{n}}
 X_{\mathbf{k}}^{''}$. The events:
\begin{equation*}
\begin{split}
&B^1_{\mathbf{n}}=\{ S_{\mathbf{n}}^{''}<-\varepsilon a_{\mathbf{n}}\}, \enskip
  B^2_{\mathbf{n}}=\{ \text{at least two} \enskip   \mathbf{k}, \enskip \mathbf{k}\leq\mathbf{n}:  -c_{\mathbf{n}}<
  X_{\mathbf{k}}<-b_{\mathbf{n}}\}\\&
  B^3_{\mathbf{n}}=\{ \text{at least one} \enskip   \mathbf{k}, \enskip \mathbf{k}\leq\mathbf{n}:
  X_{\mathbf{k}} \leq -c_{\mathbf{n}}\}. \enskip \text{and}  \enskip B_{\mathbf{n}}=\{ S_{\mathbf{n}}<
  -(\varepsilon +2\delta) a_{\mathbf{n}}\}
  \end{split}
  \end{equation*}
  satisfy the inclusion $B_{\mathbf{n}} \subset B^1_{\mathbf{n}} \cup A^2_{\mathbf{n}} \cup B^3_{\mathbf{n}}.$\\
  Likewise as the proof of (\ref{prA1}), (\ref{prA2}), (\ref{prA3}), we can show that
   \begin{equation} \label{prBn}
\sum  \limits _{ \mathbf{ n}} \mid  \mathbf{ n} \mid ^{(r/2-2)}P(B_{\mathbf{n}})  < \infty
 \end{equation}
\noindent
 By (\ref{prA1}), (\ref{prA2}), (\ref{prA3}) and (\ref{prBn}), eventually we have
$$\sum  \limits _{ \mathbf{ n}} \mid  \mathbf{ n} \mid ^{(r/2-2)} P(\mid S_{\mathbf{n}} \mid >
  (\varepsilon +2\delta) a_{\mathbf{n}})  < \infty\quad \text{for}\quad
  \varepsilon > \sigma_{1}\frac{1+\delta}{1-\delta}\sqrt{r-2} \quad \text{and all }\quad \delta>0.$$
Arbitrariness of $\delta$, allows us to conclude the implication  (\ref{momentlog} ) $\Rightarrow$ (\ref{hsulog}).
The implication (\ref{hsulogmax})$\Rightarrow$ (\ref{momentlog2}) one can prove similarly as the implication
(\ref{hsumax})$\Rightarrow$ (\ref{moment}).
 \end{pr}\\

 \medskip
 \noindent
At the end of this Section, we present one more result, which is
an extension  of some results of  Sung \cite {Sung} and  Dehua et al \cite{dehua}, to ND random fields.\\
Suppose, that $\{  \mathbf{k}_{\mathbf{n}}, \mathbf{n} \in \mathbb{N}^d \}$ is a family of lattice points of $ \mathbb{N}^d.$
\begin{theorem}
Let $\{X_{\mathbf{n},\mathbf{i}},  \mathbf{i} \leq  \mathbf{k}_{\mathbf{n}},  \mathbf{n} \in \mathbb{N}^d \}$
be an array of rowwise ND random variables with $EX_{\mathbf{n},\mathbf{i}}=0$ and
$E\mid X_{\mathbf{n},\mathbf{i}}\mid ^r < \infty$ for $1 \leq r \leq 2$,
$ \mathbf{i} \leq  \mathbf{k}_{\mathbf{n}}$ and  $\mathbf{n} \in \mathbb{N}^d.$  Furthermore assume,
that $\{a_\mathbf{n}, \mathbf{n} \in \mathbb{N}^d\}$
is a sequence of nonnegative constants. If the following conditions hold:
\begin{itemize}
  \item $\sum\limits _{\mathbf{n}}a_{\mathbf{n}}\sum \limits_{\mathbf{i} \leq  \mathbf{k}_\mathbf{n}}
  P(\mid X_{\mathbf{n},\mathbf{i}}\mid> \epsilon) < \infty $ for all $\epsilon>0$
  \item there exist $j > 0$ such that
 $
  \sum \limits_{\mathbf{n}}a_{\mathbf{n}}\left(\sum \limits_{\mathbf{i} \leq  \mathbf{k}_\mathbf{n}}
  E\mid X_{\mathbf{n},\mathbf{i}}\mid ^r  \right)^j < \infty,
 $
\end{itemize}
then $$ \sum \limits_{\mathbf{n}}a_\mathbf{n}P(\mid \sum \limits_{\mathbf{i} \leq  \mathbf{k}_\mathbf{n}}
X_{\mathbf{n},\mathbf{i}}\mid > \epsilon) < \infty \quad \text{ for all} \quad  \epsilon >0.$$
\end{theorem}
\begin{pr}By straightforward application of Lemma \ref{hofjond}.
\end{pr}

\section{Martingale random fields}
In introduction, we have given fundamental definition of martingale random field. It is known, that  we can't obtain any
sensible results for multi-parameter
martingale without any additional conditions for the  filtration. This brings us  to commutation hypothesis - also
 known as (F4) and some others:
 \begin{description}
  \item[(F3)] $\mathfrak{F}_0 $ contains all zero events of $\mathfrak{F},$
   \item[(F4)]   $\bigwedge_{ \mf{k}, \mf{n} \in \nd}$ and any bounded, $\mathfrak{F}_{\mf{k}}$-measurable random variable Y \\
   $E(Y \mid \mathfrak{F}_{\mf{n}})=E(Y \mid \mathfrak{F}_{\mf{n}\wedge \mf{k}})$ a.s.
 \end{description}
 The following notions help us to recall  definition of  strong martingale random field and j-martingale, which we exploit in
  this Section.\\
Let $J\subseteq \{1,2,...,d\}$,  $CJ= \{1,2,...,d\}\setminus J$ and
for a given $(n_1,...,n_d) \in N^d,$  denote $\mathfrak{F}_{{\mathbf{n}}}^{J}=\bigvee_{\left(n_j \in N, j
\in CJ \right)} \mathfrak{F}_{n}$.
 For example,  $ \mathfrak{F}_\mathbf{n}^{1,2}=\bigvee_{n_3,...,n_d \in N} \mathfrak{F}_{\mathbf{n}}$.\\
Thus, we can have equivalent form of (F4) condition (cf. Corollary 1of \cite{bor} ):
\begin{description}
  \item[(F4')]for any bounded random variable Y
  $$\bigwedge \limits _{\mathbf{n} \in N^d} \bigwedge \limits _{J \subseteq \{1,2,...,d\}}
   \quad  E(Y \mid \mathfrak{F}_\mathbf{n}^J \mid \mathfrak{F}_\mathbf{n}^{C\mathbf{J}})=E( Y \mid \mathfrak{F}_{\mathbf{n}}).$$
\end{description}
 Let us put $\mathcal{G}_{\mathbf{n}}=\bigvee \limits_{j=1}^d\mathfrak{F}_{\mathbf{n}}^j$ and
   $\widetilde{\mathfrak{F}}_{\mathbf{n}-\mathbf{1}}=
   \mathcal{G}_{\mathbf{n}-\mathbf{1}} \wedge \mathfrak{F}_{\mathbf{n}}$,
    where $\mathbf{n}-\mathbf{1}=(n_1-1,n_2-1,...,n_d-1)$\\
Furthermore, we need  the following conditions:
\begin{description}
  \item[(X1)] $X_{\mathbf{n}}=0$  if $ \mid \mathbf{n}\mid =0,$
   \item[(X2)] the family of random variables   $\{ X_{\mf{n}}, \in \nd\}$ is measurable with respect to family of
    $\sigma$- algebras $\{\mathfrak{F}_{\mathbf{n}}, \mathbf{n} \in \nd\},$
  \item[(X3)]  $E(X_{\mathbf{n}}\mid \widetilde{\mathcal{{F}}}_{\mathbf{n}-\mathbf{1}}) \leq 0$  \enspace
  for very $\mathbf{n} \in N^d,$
  \item[(X3')]  $E(X_{\mathbf{n}}\mid \widetilde{\mathcal{{F}}}_{\mathbf{n}-\mathbf{1}}) = 0$\enspace
  for very $\mathbf{n} \in N^d,$
   \item[(F5)] $E(Y \mid \mathfrak{F}_{\mathbf{n}}\mid \mathcal{G}_{\mathbf{n}-\mathbf{1}})=
   E(Y \mid \widetilde{\mathfrak{F}}_{\mathbf{n}-\mathbf{1}})$
   \enspace    for every $\mathbf{n} \in N^d$ and any bounded random variables Y.
\end{description}
An integrable  family of random variables $\{(X_{\mf{n}}, \mathfrak{F}_{\mf{n}}), \mf{n} \in \nd \}$
satisfying condition
(X2) is:
  \begin{itemize}
    \item  strong martingale differences iff   $E(X_{\mathbf{n}}\mid \mathcal{G}_{\mathbf{n}-\mathbf{1}})=0$ a.s.,
    \item j-martingale differences iff $(X_{\mf{n}},  \mathfrak{F}^j_{\mf{n}})$ is a one parameter martingale
    differences
    with respect to coordinate $n_j$.
  \end{itemize}
Fuk-Nagaev inequality for martingale random fields was proved by Lagodowski \cite{lag} in the case $d=2$ and
extended to the
case $d \geq 2$ by Borodhikin \cite {bor}, both authors have obtained theorems for the bounded  second
conditional moments.
We complete these results to the arbitrary r-th  conditional absolute moment, $1 \leq r \leq 2.$\\
Let us assume, that there exist fields of positive numbers $\{b_{\mathbf{k}}^r, \enskip \mathbf{k} \in \mathbb{N}^d\}$,
 $\{d_{\mathbf{k}} \enskip \mathbf{k} \in \mathbb{N}^d\}$, $\{\wt{\lambda}_{\mathbf{k}},
 \enskip \mathbf{k} \in \mathbb{N}^d\}$
 and $\{\wt{m}_{\mathbf{k}}^r, \enskip \mathbf{k} \in \mathbb{N}^d\}$ such that
\begin{equation} \label{ogrlambda}
\begin{aligned}
&E(\mid X_{\mathbf{n}}\mid^rI[\mid X_{\mathbf{n}} \mid \leq  y]\mid\widetilde{\mathfrak{F}}_{\mathbf{n}-\mathbf{1}})
\leq
b_{\mathbf{k}}^r,
&E (X_{\mathbf{k}}I[X_{\mathbf{k}}>- y]\mid \widetilde{\mathfrak{F}}_{\mathbf{n}-\mathbf{1}}) \leq
d_{\mathbf{k}}, \\
&E(\mid X_{\mathbf{n}}\mid ^r \mid  \widetilde{\mathfrak{F}}_{\mathbf{n}-\mathbf{1}}) \leq \wt{m}_{\mathbf{k}}^r,
&E (\mid  X_{\mathbf{k}} \mid I[\mid X_{\mathbf{k}}\mid > y]\mid
\widetilde{\mathfrak{F}}_{\mathbf{n}-\mathbf{1}}) \leq \wt{\lambda}_{\mathbf{k}}
\end{aligned}
\end{equation}
for every ${\mathbf{k} \in \mathbb{N}^d}$ and denote\\
$B_{\mathbf{k}}^r=\sum\limits_{\mathbf{k} \leq \mathbf{n}}b_{\mathbf{k}}^r$,
$  D_{\mathbf{n}}=\sum \limits _{\mathbf{k} \leq \mathbf{n}}d_{\mathbf{k}}$,
 $\wt{M}_{\mathbf{n}}^r=\sum\limits_{\mathbf{k} \leq \mathbf{n}}\wt{m}_{\mathbf{k}}^r$ and
  $\wt{\Lambda}_{\mathbf{k}}=\sum\limits_{\mathbf{k} \leq \mathbf{n}}\wt{\lambda}_{\mathbf{k}}$\\

\begin{theorem} \label{fukma}
Suppose, that family of $\sigma$- algebras $\{\mathfrak{F}_{\mathbf{n}}, \mathbf{n} \in \nd\}$  satisfies
condition (F1),(F3),(F4)
and (F5) in  the case  $d>2$, family of random variables $\{ X_{\mf{n}}, \in \nd\}$ satisfies
 conditions  (X1)-(X3), (\ref{ogrlambda}) and let $x,y >0$, $1 \leq r \leq 2,$ \\
then the following inequalitie holds
\begin{equation} \label{fukma1}
\begin{aligned}
&P(\max \limits_{\mathbf{k}\leq \mathbf{n}}S_{\mathbf{k}} \geq x) \leq
   P(\ma X_{\mathbf{k}} \geq y)\\
&+ e^{d-1}\exp\left\{\frac{x}{y}-\left(\frac{x-D_{\mathbf{n}}}{y}
    +\frac{B_{\mathbf{n}}^r}{y^r}\right)\ln\left[\frac{xy^{r-1}}{B_{\mathbf{n}}^r}+1\right]\right\},
\end{aligned}
\end{equation}
 if we assume (X3') instead of (X3), we have

\begin{equation} \label{fukma2}
\begin{aligned}
&P(\max \limits_{\mathbf{k}\leq \mathbf{n}}\mid S_{\mathbf{k}} \mid  \geq x) \leq
   P(\ma \mid X_{\mathbf{k}}\mid \geq y)\\
&+ 2e^{d-1}\exp\left\{\frac{x}{y}-\left(\frac{x-\wt{\Lambda}_{\mathbf{n}}}{y}
    +\frac{\wt{M}_{\mathbf{n}}^r}{y^r}\right)\ln\left[\frac{xy^{r-1}}{\wt{M}_{\mathbf{n}}^r}+1\right]\right\}.
\end{aligned}
\end{equation}

    \end{theorem}
\begin{pr}(sketch)  Let us put
 \begin{equation*} \begin{split}
 & \widetilde{X_{\mathbf{k}}}=X_{\mathbf{k}}I[X_{\mathbf{k}} \leq y], \enspace \widetilde{S_{\mathbf{k}} }=
  \sum \limits _{\mathbf{k} \leq \mathbf{n}}\widetilde{X_{\mathbf{k}}}\enspace \text{and}\\&
  Z_{\mathbf{k}}=\widetilde{X_{\mathbf{k}}}-E(\widetilde{X_{\mathbf{k}}}
  \mid \widetilde{\mathfrak{F}}_{\mathbf{n}-\mathbf{1}}),
  \enspace  T_{\mathbf{n}}=\sum \limits _{{\mathbf{k} \leq \mathbf{n}}}Z_{\mathbf{k}}.
  \end{split}
\end{equation*}
Obviously, we have
 \begin{equation} \label{nier1}
P(\max \limits_{\mathbf{k} \leq \mathbf{n}}S_{\mathbf{k}}\geq x) \leq P(\max \limits_{\mathbf{k}
 \leq \mathbf{n}}\widetilde{S_{\mathbf{k}}}\geq x)+
P(\max \limits_{\mathbf{k} \leq \mathbf{n}}X_{\mathbf{k}}\geq y).
\end{equation}
From (X3) implies, that $Z_{\mathbf{k}} \geq \widetilde{X_{\mathbf{k}}}$  a.s. and since $\alpha >1$, $h>0$

\begin{equation}  \label{nier2}
P(\max \limits_{\mathbf{k} \leq \mathbf{n}}\widetilde{S_{\mathbf{k}}}\geq x) \leq P(\max \limits_{\mathbf{k}
\leq \mathbf{n}}e^{\alpha h T_{\mathbf{k}}}\geq e^{\alpha h x}).
\end{equation}
\\
Let us observe, that $\{ (e^{\alpha h T_{\mathbf{k}}}, \widetilde{\mathfrak{F}}_{\mathbf{k}}),
\mathbf{k} \leq \mathbf{n} \}$ is positive submartingale.\\
Denote $\mathbf{k}(j)=(k_1,k_2,\cdots,k_{j-1},k_{j+1},\cdots,k_d)$ for $\mathbf{k} \in \mathbb{N}^d$
and $1 \leq j \leq d$, thus\\
$\{\max \limits_{\mathbf{k}(d) \leq \mathbf{n}(d)} e^{\alpha h T_{\mathbf{k}}}, 1 \leq k_d \leq n_d\}$
 is positive d-sumbartingale with respect to
$\{ \mathfrak{F}_{\mathbf{k}}^{d}, \mathbf{k} \leq \mathbf{n}\}.$
By application of standard Doob inequality to d-submartingale and Doob inequality for submartingale random field,
cf.Shorack et al.  \cite{sho}
\begin{equation}  \label{nier3}
\begin{split}
&P(\max \limits_{\mathbf{k} \leq \mathbf{n}}e^{\alpha h T_{\mathbf{k}}}\geq e^{\alpha h x}) \leq e^{-\alpha h x}
E\left(\max \limits_{\mathbf{k}(d) \leq \mathbf{n}(d)}( e^{ h T_{\mathbf{k}(d)n_d}})^{\alpha} \right)\\&
\leq \left( \frac{\alpha}{\alpha-1} \right)^{\alpha(d-1)}e^{-\alpha h x}Ee^{\alpha h T_{\mathbf{n}}}.
\end{split}
\end{equation}
Furthermore, we need estimations:
\begin{itemize}
  \item $
\begin{aligned}[t]
&E(e^{\alpha h \widetilde{X_{\mathbf{k}}}}\mid \widetilde{\mathfrak{F}}_{\mathbf{k}-\mathbf{1}})=
E(e^{\alpha h \widetilde{X_{\mathbf{k}}}} I[\widetilde{X_{\mathbf{k}}}< -y]\mid
\widetilde{\mathfrak{F}}_{\mathbf{k}-\mathbf{1}})+\\&
E(e^{\alpha h \widetilde{X_{\mathbf{k}}}} I[\mid\widetilde{X_{\mathbf{k}}} \mid \leq y]\mid
\widetilde{\mathfrak{F}}_{\mathbf{k}-\mathbf{1}})=I_{10}+I_{11},
\end{aligned}
$
  \item$I_{10} \leq E( I[X_{\mathbf{k}}< -y]\mid \widetilde{\mathfrak{F}}_{\mathbf{k}-\mathbf{1}}),$
  \item  $
\begin{aligned}[t]
&I_{11} \leq  \frac{e^{\alpha h y - 1 - \alpha h y}}{y^2} E\left(X_{\mathbf{k}}^2I[ 0<\mid X_{\mathbf{k}}  \mid \leq y] \mid
\widetilde{\mathfrak{F}}_{\mathbf{k}-\mathbf{1}}\right)+ \alpha h E( X_{\mathbf{k}}I[\mid X_{\mathbf{k}}  \mid \leq y] \mid
\widetilde{\mathfrak{F}}_{\mathbf{k}-\mathbf{1}})+\\&
E( I[\mid X_{\mathbf{k}}  \mid \leq y] \mid \widetilde{\mathfrak{F}}_{\mathbf{k}-\mathbf{1}}) \leq
\frac{e^{\alpha h y - 1 - \alpha h y}}{y^r}E \left(\mid X_{\mathbf{k}} \mid ^rI[ 0<\mid X_{\mathbf{k}} \mid \leq y] \mid
\widetilde{\mathfrak{F}}_{\mathbf{k}-\mathbf{1}}\right)+\\&
\alpha h E( X_{\mathbf{k}}I[\mid X_{\mathbf{k}}  \mid \leq y] \mid
\widetilde{\mathfrak{F}}_{\mathbf{k}-\mathbf{1}}).
\end{aligned}
$
\end{itemize}
Thus
 \begin{equation}
\begin{split}
&E(e^{\alpha h Z_{\mathbf{k}}}\mid \widetilde{\mathfrak{F}}_{\mathbf{k}-\mathbf{1}})\leq \\&
e^{-\alpha h  E( X_{\mathbf{k}}I[X_{\mathbf{k}} \leq y]
\mid\widetilde{\mathfrak{F}}_{\mathbf{k}-\mathbf{1}})}\exp\Big\{ \frac{e^{\alpha h y - 1 - \alpha h y}}{y^r}E \left(\mid X_{\mathbf{k}} \mid ^rI[ 0<\mid X_{\mathbf{k}} \mid \leq y] \mid
\widetilde{\mathfrak{F}}_{\mathbf{k}-\mathbf{1}}\right) +\\&
 \alpha h E( X_{\mathbf{k}}I[\mid X_{\mathbf{k}}  \mid \leq y] \mid
\widetilde{\mathfrak{F}}_{\mathbf{k}-\mathbf{1}})\Big\} \leq \exp\Big\{ \frac{e^{\alpha h y - 1 - \alpha h y}}{y^r}
E \left(\mid X_{\mathbf{k}} \mid ^rI[ 0<\mid X_{\mathbf{k}} \mid \leq y] \mid
\widetilde{\mathfrak{F}}_{\mathbf{k}-\mathbf{1}}\right) -\\&
 \alpha h E( X_{\mathbf{k}}I[ X_{\mathbf{k}}  < -y] \mid
\widetilde{\mathfrak{F}}_{\mathbf{k}-\mathbf{1}})\Big\} \leq \exp\Big\{ \frac{e^{\alpha h y - 1 - \alpha h y}}{y^r}b_{\mathbf{k}}^r+
 \alpha h d_{\mathbf{k}}\Big\}.
\end{split}
\end{equation}
Now, furnishing $\{\mf{k}:\mf{k} \leq \mf{n}\}$ with a  total order and using property (F5), we have
 \begin{equation} \label{nier4}
 Ee^{\alpha h T_{\mathbf{n}}} \leq  \exp \Big\{ \frac{e^{\alpha h y - 1 - \alpha h y}}{y^r}
 \sum \limits _{\mathbf{k} \leq \mathbf{n}} b_{\mathbf{k}}^r+ \alpha h
  \sum \limits _{\mathbf{k} \leq \mathbf{n}} d_{\mathbf{k}}\Big\}.
 \end{equation}
Combining (\ref{nier1}), (\ref{nier2}), (\ref{nier3}) and (\ref{nier4}) we get
\begin{equation}\begin{split}
&P(\max \limits_{\mathbf{k}\leq \mathbf{n}}S_{\mathbf{k}} \geq x) \leq \\&
P(\max \limits_{\mathbf{k} \leq \mathbf{n}}X_{\mathbf{k}}\geq y)    +  \left( \frac{\alpha}{\alpha-1} \right)^{\alpha(d-1)}
e^{-\alpha h x}\exp \Big\{ \frac{e^{\alpha h y - 1 - \alpha h y}}{y^r}B_{\mathbf{k}}^r+
 \alpha h D_{\mathbf{k}} \Big\}    \leq \\&
 P(\max \limits_{\mathbf{k} \leq \mathbf{n}}X_{\mathbf{k}}\geq y)    +  e^{d-1}
\exp \Big\{ \frac{e^{\alpha h y - 1 - \alpha h y}}{y^r}B_{\mathbf{k}}^r+
 \alpha h D_{\mathbf{k}}-\alpha h x\Big\}.
 \end{split}
\end{equation}
Setting,  $\alpha h = \frac{1}{y}\ln[\frac{xy^{r-1}}{B^r_{\mf{n}}1}+1]$ one can obtain (\ref{fukma1}).\\
To prove (\ref{fukma2}), we set:
$Y_{\mf{k}}=-\x$ and $U_{\mf{n}}=\su Y_{\mf{k}}$. Obviously, $\{(U_{\mf{n}}, \mathfrak{F}_{\mf{n}}), \mf{n}\in \nd \}$
is martingale random field satisfying assumption of our theorem. Furthermore, denote
$\wt{Y}_{\mf{k}}=Y_{\mf{k}}I[Y_{\mf{k}} \leq y]$ and $\wt{U}_{\mf{n}}=\su \wt{Y}_{\mf{k}}$.
By standard estimation we have
\begin{equation*}
P(\max \limits_{\mathbf{k} \leq \mathbf{n}}\mid S_{\mathbf{k}} \mid \geq x) \leq  P(\max \limits_{\mathbf{k} \leq \mathbf{n}}\mid X_{\mathbf{k}}\mid \geq y) +
P(\max \limits_{\mathbf{k}\leq \mathbf{n}}\widetilde{S_{\mathbf{k}}}\geq x)+P(\max \limits_{\mathbf{k}
 \leq \mathbf{n}}\widetilde{U_{\mathbf{k}}}\geq x),
\end{equation*}
then  similarly, as in the first part of the proof we obtain  (\ref{fukma2}).

\end{pr}
\begin{lemma}\label{hofjoma}
Let $\{(X_{\mf{n}}, \mathfrak{F}_{\mf{n}}), \mf{n} \in \nd \}$ satisfies assumption of Theorem  \ref{fukma},
then there exist constant $C>0$ such that  for every $x>0$ and $j>0$
\begin{equation*}
P(\max \limits_{\mathbf{k}\leq \mathbf{n}}\mid S_{\mathbf{k}} \mid  \geq x) \leq
   P(\ma \mid X_{\mathbf{k}}\mid \geq \frac{x}{j})+C\left( \frac{1}{x^r}\wt{M}_{\mathbf{n}}^r\right)^j.
\end{equation*}
\end{lemma}
\noindent
The proof of this lemma is similar to those of Lemma \ref{hofjond}, thus we omit it.\\
Application of  Lemma \ref{hofjoma} gives the following two theorems. The first one is martingale random field
version of Theorem 1.3 ((1.10) $\Rightarrow$(1.11)) of Gut et al. \cite{GutHsu} and the  latter, an extension of
Theorem 1(b) and Theorem 2 of Ghosal and Chandra  \cite {gho} to  martingale random field with weaker moment restriction.
\\
\begin{theorem}
Let $\{(X_{\mf{n}}, \mathfrak{F}_{\mf{n}}), \mf{n} \in \nd \}$ satisfies assumption of Theorem  \ref{fukma}
and WMD condition,  moreover
suppose that $\alpha_1> \frac{1}{2}$,  $\alpha_1r > 1$ and
there exist constant M depend only on r and $\mf{n}_{\mf{1}} \in \nd $ such that $ \frac{1}
{\mid \mf{n} \mid}\wt{M}_{\mathbf{n}}^r \leq M$ for any $\mf{n}\geq \mf{n}_{\mf{1}} $ , thus if
\begin{equation} \label{moment2}
E\mid\xi\mid^r (\log_{+}\mid \xi\mid)\mid^{p-1}<\infty,
\end{equation}
then
\begin{equation} \label{hsu2}
\sum\limits_{\mathbf{n}}\mid \mathbf{n}\mid^{\alpha_1r-2}P(\max \limits_{\mf{k}\leq \mf{n}}\mid S_{\mathbf{k}}\mid>\mid \mathbf{n}^{\boldsymbol{\alpha}}\mid \varepsilon)< \infty   \quad \text{for all}\quad  \varepsilon >0.
\end{equation}

\end{theorem}
\begin{pr}
Likewise the proof  of (\ref{moment}) $\Rightarrow$(\ref{hsu}).
\end{pr}\\

\noindent
Now, suppose that $\{  \mathbf{k}_{\mathbf{n}}, \mathbf{n} \in \mathbb{N}^d \}$
is a family of lattice points of $ \mathbb{N}^d$.
\begin{theorem}
Let $\{X_{\mathbf{n},\mathbf{i}},  \mathbf{i} \leq  \mathbf{k}_{\mathbf{n}},  \mathbf{n} \in \mathbb{N}^d \}$
 be a d-dimensional array of rowwise martingale differences with respect to family of $\sigma$-algebras\\
 $\{\mathfrak{F}_{\mathbf{n},\mathbf{i}},  \mathbf{i} \leq  \mathbf{k}_{\mathbf{n}},  \mathbf{n} \in \mathbb{N}^d \}$
 satisfying assumptions of Theorem  \ref{fukma}.
  Furthermore, there exist families of non-negative constants $\{a_\mathbf{n}, \mathbf{n} \in \mathbb{N}^d\}$ and
  $\{\widetilde{M_{ \mathbf{k}_\mathbf{n}}^r }, \mf{n} \in \nd\}
  \}$ such that
  $\sum \limits _{\mathbf{j} \leq \mathbf{k}_\mathbf{n}}
  E(\mid X_{\mathbf{n},\mf{j}}\mid^r \mid \widetilde{\mathfrak{F}}_{\mf{n},\mathbf{j}-\mathbf{1}}) \leq
  \widetilde{M_{ \mathbf{k}_\mathbf{n}}^r }.$
  If the following conditions hold:
\begin{itemize}
  \item $\sum\limits _{\mathbf{n}}a_{\mathbf{n}}   P(\max \limits_{\mathbf{i} \leq  \mathbf{k}_\mathbf{n}}
  \mid X_{\mathbf{n},\mathbf{i}}\mid > \epsilon) < \infty $ for all $\epsilon>0$
  \item there exist $j > 0$ such that
 $ \sum \limits_{\mathbf{n}}a_{\mathbf{n}}\left(\widetilde{M_{ \mathbf{k}_\mathbf{n}}^r } \right)^j < \infty,$
\end{itemize}
then  \begin{equation*}
 \sum \limits_{\mathbf{n}}a_\mathbf{n}P(\max \limits_{\mathbf{l} \leq  \mathbf{k}_\mathbf{n}}
 \mid \sum \limits_{\mathbf{i} \leq
 \mathbf{l}} X_{\mathbf{n},\mathbf{i}} \mid > \epsilon) < \infty \quad \text{ for all }\quad \epsilon >0.
  \end{equation*}
\end{theorem}
\medskip
\noindent
{\bf Negatively associated random fields - some comments.}\\
Fuk-Nagaev inequality for sequences of negatively associated random variables can be proved by application of the
comparison theorem which has been obtained by Shao \cite{shao}.
\begin{theorem}
Let $\{X_i, 1\leq i \leq n \}$ be a negatively associated sequence and let $\{X_i^*, 1 \leq i \leq n \}$ be a sequence of independent random variables such that $X_i^*$ and $X_i$ have the same distribution for each $i=1,2,\ldots,n$. Then
\begin{equation}
Ef\left(\max \limits_{1 \leq k \leq n } \sum \limits_ {i=1}^kX_i \right) \leq Ef\left(\max \limits_{1 \leq k \leq n }\sum \limits _{i=1}^k X_i^*\right)
\end{equation}
for any convex and non-decreasing function f on $\mathbf{R}^1$, whenever the expectation on the right side exist.
\end{theorem}
In the case $d\geq 2$, Bulinski and Suquet  (cf, Theorem 2.12 of \cite{bul}) have proved, that the comparison
theorem does not hold in general for maximum of sums of NA random field.
\begin{theorem}
Let $ f : \mathbf{R}\rightarrow \mathbf{R}$ be a function such that $f(1) > f(0)$ (in
particular, strictly increasing). Then, for any $d > 1$ there exists a NA random field
 $X = \{X_{\mathbf{j}} ; \mathbf{j}\in  Z^d\}$ and a multiindex  $\mathbf{n} \in N^d $ such that
\begin{equation}
Ef\left(\max \limits_{ \mathbf{k} \leq \mathbf{n} } \sum \limits_ {  \mathbf{i} \leq \mathbf{k} }
X_{\mathbf{i}} \right)> Ef\left(\max \limits_{ \mathbf{k} \leq \mathbf{n} }
\sum \limits _{  \mathbf{i} \leq \mathbf{k} }X_{\mathbf{i}}^*\right),
\end{equation}
where $X^{*} = \{X^{*}_{\mathbf{j}} ; \mathbf{j}\in  Z^d\}$ is decuopled version of $X$.
\end{theorem}
Furthermore, Shao has used in his  proof,  the  maximal inequality for non-negative supermartingale
which is not true for supermartingale random
fields. Thus, the maximal Fuk-Nagaev inequality for NA random fields is open question.

\end{document}